\documentclass[twoside,12pt]{article}
\usepackage{amsmath,amssymb,epsfig, amsthm}
\date{}
\newcommand{\CL}{{\cal L}}

\def\phi{\varphi}
\def\epsilon{\varepsilon}

\newcommand{\BR}{\mathbb R}
\newcommand{\BC}{\mathbb C}
\newcommand{\BH}{\mathbb H}
\newcommand{\BO}{\mathbb O}
\newcommand{\BF}{\mathbb F}
\newcommand{\BD}{\mathbb D}

\newcommand{\bthm}{\begin{Theorem}}
\newcommand{\ethm}{\end{Theorem}}
\newcommand{\ble}{\begin{Lemma}}
\newcommand{\ele}{\end{Lemma}}
\newcommand{\bprop}{\begin{Proposition}}
\newcommand{\eprop}{\end{Proposition}}
\newcommand{\bcor}{\begin{Corollary}}
\newcommand{\ecor}{\end{Corollary}}

\begin{document}
\title{Helmut Salzmann and his  legacy}
\author{Rainer L\"owen}
\maketitle
\thispagestyle{empty}

\begin{abstract}
We describe the development of the mathematics of Helmut R. Salzmann 
(3. 11. 1930 -- 8. 3. 2022) and the main difficulties he was facing, 
documenting his lifelong productivity and his far reaching influence.
We include a comprehensive bibliography of his work.


\end{abstract}

Helmut Salzmann left an enormous heritage. It is impossible to do justice to his 
oeuvre within a few pages. Instead of attempting this, we shall focus here on 
the development of his work and the main difficulties he was facing, highlighting a few 
points of special interest.

\section {Previous work that inspired Salzmann}

The first complete axiomatic description of the affine plane over the reals 
together with its group of motions (or isometries) was given by David Hilbert in his 
seminal book `Grundlagen der Geometrie' \cite{hilb} (see also \cite{fr57b} 
for further historical information). 
This was an axiom system
in the spirit of Euclid, with axioms of three types: axioms of incidence, 
yielding an affine plane, axioms of continuity, and axioms of congruence. 
The axioms of continuity were expressed in terms of order structures. One axiom
stipulated that the order on lines should be Archimedean, another one asked 
for maximality of the models. The existence of a metric was \it not \rm
presupposed, instead the axioms of congruence ensure the existence of a group of 
automorphisms (also called collineations) that is not too small 
(e.g., transitive on flags, that is, incident 
point-line pairs), yet  not too big (e.g., the stabilizer of a flag has oder 2).
Hilbert gave several examples of geometries violating some few of his axioms. 
Among them are the Moulton planes, $M_k$, $k>1$,
obtained from the real plane by changing all 
lines of negative slope $s$ in the right half plane such that the right part has slope $ks$ and 
such that the line remains connected. All of Hilbert's examples have point set $P=  \BR^2$, 
and their lines are closed subsets homeomorphic to $ \BR$.

\

Hilbert's book together with the examples like Moulton's triggered intense 
work on affine and projective planes. 
It emerged that after choosing a quadrangle as frame of reference, affine 
(inhomogeneous) coordinates 
can be introduced and be equipped with an 
algebraic structure, 
but in general, this is not a (skew) field. Instead, one gets a single ternary 
operation $\tau$ on an affine line $K$. The affine plane then appears as the point 
set $K\times K$ with lines of the form
  $$L_{s,t} = \{(x,y) \vert \ y = \tau(s,x,t) \},$$
where $s,t\in K$.  In the presence of 
some geometric regularity this operation can be split up into addition and 
multiplication as $\tau (s,x,t) = sx + t$ (such ternaries are called linear), and further 
regularity is required for every single algebraic law (like associativity etc.) 
to hold. Geometric
regularity here means either closure theorems like that of Desargues or special 
cases thereof or, 
equivalently, existence of 
groups of central collineations (fixing all lines through some point) that 
are as transitive as possible.
The unrivalled, comprehensive book \cite{pick} by Pickert on these matters had just 
been completed, and Salzmann had the opportunity to do the proofreading and thus to absorb
its full contents, see also \cite{sz16a}. The book appeared in 1955.

\

Surely the first time that some vague notion of topological projective geometries appeared
was in in Kolmogoroffs article \cite{kol} of 1932.
Pontryagin \cite{pont} had just shown that there are only three locally compact, 
connected topological 
skew fields, namely the real numbers $\BR$, the complex numbers $\BC$, and
the quaternions $\BH$. 
Kolmogoroff remarked that this implies a characterization 
of finite dimensional projective spaces over these skew fields. In a vague fashion, 
he assumed continuity of geometric operations and suggested that this corresponds to continuity
properties of associated coordinate skew fields.

\section {The beginnings}

One of Salzmann's two 1955 papers \cite{sz-z} was joint work with Zeller 
on applications of Baire category to existence of
continuous, nowhere differentiable functions. The other one \cite{sz55} 
started his lifelong devotion 
to topological projective planes, i.e., projective planes with continuous geometric operations. 
This approach was designed to replace Hilbert's continuity axioms. The paper mentioned undertook
a very thorough investigation of the topological properties, in particular, 
connectedness, assuming nothing except that the topology should not be the coarsest. 
The results were incorporated in Pickert's book.

\

Around the same time, two other authors had started investigations of topological 
projective planes. In 
Russia, Skornjakov wrote an article \cite{sk54} on the subject. Salzmann became 
aware of this while proofreading his own second paper. Next, there is 
Freudenthal's article \cite{fr57a},
which contains a visionary conjecture: lines of a projective plane 
should be spheres. Together with a (slightly younger) famous result by Adams on spheres with a 
continuous multiplication \cite{adams}, this would imply that lines are spheres 
of dimension $l \in \{1,2,4,8\}$, like in the 
case of the \it classical planes \rm over $\BR, \BC, \BH,\BO$, where $\BO$ denotes 
the non-associative division algebra of octonions. Freudenthal proved his conjecture
this under the assumption that lines are manifolds. He used the contracting effect of 
multiplication by the elements of a path joining the additive and multiplicative units $0$ 
and $1$. I think Salzmann once mentioned that he
met Freudenthal at the Amsterdam ICM in 1954. 

\

All the research at this time, including Salzmann's doctoral thesis \cite{sz57}, 
centered around the 
connections between topological properties of a plane and of its lines, and between 
the continuity properties of an affine plane and of its projective completion, as 
well as those of the coordinatizing ternary fields. All these things turned out to be much 
more subtle than suggested by Kolmogoroff, and most of the proofs required extra assumptions of 
a topological or algebraic nature. It was not just a matter of finding general proofs: 
much later, Eisele  constructed a 
series of counterexamples related to these questions, see references given in \cite{CPP}. 
In one case, Eisele used model theory to achieve his goal.

\

Luckily, all the problems mentioned become manageable if one assumes that the planes are 
locally compact and connected (as we shall do henceforth). 
Still the proofs require much attention to detail, 
but all desirable results in this direction could be proved in the subsequent years, 
with Salzmann 
taking a leading role. 

\

In his doctoral thesis \cite{sz57}, Salzmann also constructed his first examples of 
non-desarguesian 
planes homeomorphic to the real projective plane. Their lines are made up of pieces of 
ordinary lines like in the Moulton planes, see figure XXX. Their coordinate ternaries 
are so-called neofields, i.e., they have commutative and associative multiplicative 
structure and satisfy 
both distributive laws, but their addition is non-associative. The existence of 
these examples answered a question of Pickert, 
who was the first thesis advisor; the second advisor was Reinhold Baer, who had just 
returned to Frankfurt from his US exile.

The study of topological field-like algebraic structures associated with projective 
planes became a lasting major occupation for Salzmann, shared by his friend Karl Heinrich 
Hofmann, who worked with a more general notion, independent of planes. 

\

Skornjakov \cite{sk57} later studied  systems of curves in the plane $\BR^2$, 
requiring only that 
these `lines' are closed subsets homeomorphic to the real line, and that any two 
distinct points are joined by a unique line. He proved that, with respect to a 
suitable topology on the line set $\CL$, joining points by a line is a continuous operation, 
and the same holds true for the intersection operation $ \CL \times  \CL \to \BR^2$. 
See \cite{CPP},
Section 31 for an accessible presentation. Salzmann took up this thread and coined a new notion
of `topological planes', consisting of a point space together with a family $\CL$ of 
closed subsets, called lines, such that any two points are joined by a unique line, the 
join and intersection operations are continuous, and the set of intersecting line pairs is open.
The last condition distinguishes such planes from, e.g., spatial geometries. Because of 
the importance of this `stability condition', topological planes in this sense were 
later called stable planes. Open subsets of stable planes (in particular, open 
subsets of topological projective planes) are again stable planes, but there
are examples that do not embed in any topological projective plane. 
A nice example of this kind is the \it hyperbola plane \rm \cite{sz62b}, compare \cite{sz63a},
\cite{sz67b}, 5.3: in the real affine 
plane, delete all lines of negative slope. Instead, take the graph of $x \to x^{-1}$ 
for $x>0$ and all its images under translation. FIGUR aus  XXX
In the comprehensive survey article \cite{sz67b}, which also covered 
the topics described in the next section, Salzmann started a systematic study of stable 
planes in the same spirit as for projective planes. 

\section {Exploiting symmetry}

Symmetries of planes had been playing a major role even up to this time 
(we mentioned central collineations), but a new era started when Salzmann began to study 
the group $\Sigma$ of all automorphisms of a 2-dimensional projective plane systematically. 
Arens \cite{arens} had  shown that groups of homeomorphisms of nice spaces become topological 
transformation groups when equipped with the compact open topology, and Montgomery and 
Zippin had completed their book \cite{mozip} on transformation groups recently. Salzmann proved 
that the stabilizer of a quadrangle 
in the group $\Sigma$ is trivial (a property later called stiffness), and this enabled him to 
embed the group in the space of all quadrangles, in order to prove its local compactness. 
Then Montgomery and Zippin came in to show that the group is a Lie group, thus opening 
a very rich toolbox. 

\

If the plane is homogeneous, i.e., if $\Sigma$ is transitive on the point set, then by 
a result of Montgomery \cite{mont} there is a compact transitive subgroup. 
Salzmann showed that this 
group must be the rotation group $\mathop{\mathrm{SO}_3}\BR$. 
This group contains many involutions 
(elements of order two), and he proved that involutions must be reflections 
(having center and axis). An isomorphism onto the real plane is then obtained by 
matching centers and axes of reflections. Incidence of points and lines is equivalent 
to commutability of the corresponding reflections, and this guarantees that the resulting 
map of planes is an isomorphism.

\

Lesser degrees of symmetry of a plane could now be measured
as the dimension of its  Lie group of automorphisms, 
and this became a leading idea for all 
future work. In a series of four papers \cite{sz62a}, \cite{sz63b}, \cite{sz64}, 
\cite{sz65}, Salzmann showed 
that for 2-dimensional planes $d = \dim \Sigma \le 8$ holds. The bound is attained for the real 
plane, and no other plane has $d\ge 5$. The planes with $d=4$ are precisely the Moulton planes
$M_k$, $k \ge 1$. The Moulton planes are something very special, even 
among \it all \rm compact connected planes. In fact, Salzmann showed in \cite{sz74} that 
no other plane shares their configuration of transitive elation groups (called Lenz type III).
The planes with $d = 3$ form three completely new families depending 
on 1, 2 and 3 parameters, respectively. One of the families consists of the \it
cartesian planes \rm depicted in Figure XXX.

\

An important step towards this result was the fact that involutions are 
reflections (as mentioned above), with the consequence that by taking products of involutions, 
elations (or in affine terminology, translations) could be generated, that is, 
central collineations with the center lying on the axis. Elations with common axis tend to 
form commutative groups, which helped to unravel the structure of $\Sigma$. 
Another powerful tool was the classification of all transitive 
groups on 1-manifolds, obtained by Brouwer as early as 1909. Here the universal covering group
$\mathop{\mathrm{SL}_2^\sim}\BR$ occurs, which acts on the Moulton planes, 
making these planes an excellent model 
of that elusive group, compare \cite{CPP}, Section 34.

\

Perhaps the nicest planes ever discovered by Salzmann are the skew-hyperbolic ones, admitting 
the standard action of the hyperbolic motion group $\mathrm{SO}_3(\BR,1)$. 
The group leaves a conic invariant, 
and the parts of lines in the interior of this conic are intervals of ordinary lines, 
while the exterior parts are orbits of the line stabilizers in the shape of hyperbola branches,
see Figure XXXX. 

\

It is worth while to review the main steps in the classification program that was 
executed here for the first time. One starts from some hypothesis on the degree of 
symmetry. Then one tries to determine the possible group structures of $\Sigma$, and, 
after that, the possible actions of those groups as transformation groups of the point 
set, paying special attention to fixed elements and stabilizers.
From this one hopes to deduce a description of the lines as subsets of the point set. 
For all possible line sets that emerge in this process, it remains finally to 
determine whether or not  they define projective planes. Thus one has to check whether or not
certain equations describing intersection points and joining lines have 
unique solutions. This last step can be extremely hard. 

\

Some of  Salzmann's articles on the classification of 2-dimensional 
projective planes were written 
during his first one-year visit to Los Angeles, where he worked with Herbert Busemann 
and with James Dugundji, who became a very close friend. With Busemann he wrote a 
paper on synthetic differential geometry \cite{bu-sz}.

\section {Dimension four}

The foundations for the study of 4-dimensional projective planes were laid in \cite{sz69b}.
Using very recent topological results, Salzmann showed first of all that there 
are no planes of dimension three, so 4-dimensional planes are next to be considered. 
They are the relatives 
of the complex projective plane and share its topological properties, 
compare \cite{CPP} 53.15. 
Much of the study of these 
planes runs along the same lines as in the 2-dimensional case, only the groups become bigger and
hence there are more possibilities to be considered. 
However, the most significant differences occur with the 
structure of ternaries and with stiffness. The field of complex numbers contains 
the real numbers as a subfield, so 2-dimensional ternaries may have 
proper closed subternaries, but they don't have to. If there is such a subternary, 
then it is homeomorphic to $\BR$ and, in particular, not discrete or 0-dimensional. 
This is obtained using Brouwer's theorem on fixed point free homeomorphisms of $\BR^2$ 
(stating in particular that they have infinite order) \cite{brouwer12}.

\

A 1-dimensional subternary corresponds to a 2-dimensional subplane, 
which is a so-called Baer subplane 
(named after Reinhold Baer). This means that every point not on the subplane is 
incident with exactly one line of the subplane. Salzmann showed by topological arguments 
that the concepts `Baer subplane' and `2-dimensional subplane' are equivalent 
in 4-dimensional planes.
The stabilizer of a quadrangle 
in the complex plane is generated by the automorphism induced by complex conjugation 
of homogeneous coordinates and thus is of order two. Note, however, that $\BC$ 
has discontinuous automorphisms, a fact that made some precautions necessary when 
dealing with the general situation.
In general, a Baer subplane of 
a 4-dimensional plane may be fixed pointwise by a nontrivial automorphism of the larger 
plane, which then has to be an involution, called a Baer involution. All involutions 
of a 4-dimensional plane are either Baer involutions, reversing the orientation of 
their fixed lines, or reflections, preserving the orientation of fixed lines. Finally, 
this leads to the fact that the stabilizer of a quadrangle has order at most two (stiffness).
A key ingredient is once more Brouwer's theorem on fixed point free homeomorphisms of $\BR^2$. 

\

In \cite{sz70}, Salzmann shows that the automorphism group $\Sigma$ of a 4-dimensional 
projective plane is locally compact. The key step is to show that a sequence of automorphisms
is uniformly convergent if it converges on the vertices of a quadrangle plus one more 
suitable point. On the basis of this fact, one obtains next that $\Sigma$ is a Lie group.
Here, Salzmann uses the result of Montgomery and Zippin \cite{mozip} that a 
locally compact group acting 
on a surface is a Lie group, applying it to the effective quotient action of a line 
stabilizer on the line (which is a 2-sphere). 
By stiffness, it follows that  $d = \dim \Sigma \le 16$ holds for the automorphism group of 
a 4-dimensional projective plane. The bound is attained for the classical complex plane.
Using results by P.A. Smith on transformation groups, Salzmann obtains that reflections 
are determined by their fixed elements. Baer involutions can be kept away because of 
their orientation behavior, hence the result that the complex plane is the only 
4-dimensional plane whose group is transitive on the point set can be proved much like in
the 2-dimensional case. 

\

In the subsequent article \cite{sz71}, it is shown that the complex plane is the only one 
satisfying $d > 8$. The proof is very involved, more so and longer than previous ones 
in Salzmann's work. The starting point is that by the Lefschetz fixed point theorem, 
every automorphism fixes a point and a line. Considering the possible dimensions $d$ 
one by one in descending order, Salzmann shows each time that a group of dimension $d$ 
fixes a point or a line. A key role is played by large commutative groups of elations and,
of course, by the structure theory of Lie groups. Actually, in the case $d = 9$, the 
proof ends by showing that the plane is a \it translation plane\rm, i.e., that 
$\Sigma$ contains a transitive group of elations with a common axis. Translation planes 
of dimension 4 were studied by D. Betten around the same time, and he showed that $d = 9$ 
is impossible \cite{betten72b}. 

\

This proof is typical of most proofs later devised by Salzmann, which become 
harder and harder over the years. One of the fundamental tools used again and again is 
the dimension formula for orbits of transformation groups; under suitable conditions, 
the topological dimension of an orbit is the dimension of the group minus that of a stabilizer. 
So large groups have large orbits or large stabilizers or both.
For Lie groups, the decomposition as a semidirect product of the solvable radical and 
a semi-simple Levi complement is a standard tool, as well as the Malcev-Iwasawa 
decomposition as a topological product of a maximal compact subgroup and a vector space.
The classification by Mostow of transitive actions on surfaces is used; in higher dimensions,
its role is taken over by Richardson's classification 
\cite{rich} of compact groups acting on the 4-sphere.

\

In 1970, both Betten \cite{betten70} and Plaumann and Strambach \cite{p-s70}
constructed proper (non-classical) translation planes by 
modifying the multiplication of $\BC$. No such thing exists in dimension two. 
In fact, Betten classified all translation planes with 
$d \ge 7$, compare \cite{CPP} Section 73. For $d = 8$, 
there are three possibilities for the group, 
leading to two single planes and a one-parameter family of planes  
\cite{betten73a}, \cite{betten73b}. 
So at this point, it was clear that $d = 8$ is the largest dimension of automorphism groups 
of non-classical planes. This did not satisfy Salzmann; he wanted to know at least 
all planes with d = 8. If the group $\Sigma$ is not simple, this was achieved in  
\cite{sz73b} by showing that the plane is a translation plane. Up to local isomorphism, 
there are three simple groups of dimension 8. The compact motion group of the complex plane
can only act point transitive, and then the plane is classical as mentioned earlier. 
The hyperbolic (non-compact) motion group of the complex plane again acts only on the complex 
plane, as was shown by Betten \cite{betten73c}. This leaves the 
automorphism group $\Sigma_\BR$ of the 
real plane as the only remaining possibility. This group does act on the complex 
plane, leaving the real (Baer-)subplane invariant. Clearly the subplane is 
completely determined by this action. In \cite{sz73c}  it is shown that in fact, the action 
of this group even determines the exterior structure, too, so that the plane is classical.
In other words, all of the most homogeneous non-classical 4-dimensional pojective planes are
translation planes which are completely known. This pattern repeats itself in 
higher dimensions, with the exception of the last case where an 
invariant Baer subplane exists: in dimensions 8 and 16, 
non-classical so-called  Hughes planes occur, see Section \ref{Hughes}.

\section {Number systems}

From what we have seen so far, it may be recognized that there is a very close connection 
between topological projective planes and number systems. In the classical cases, this 
means real and complex numbers, quaternions and octonions.  In the course of his 
research, Salzmann frequently came across questions regarding the structure of these 
number domains, and he had thus accumulated a very intimate knowledge. In oder to share this, 
he gave a two Semester lecture course entitled  `Zahlbereiche' starting in summer, 1971, and he 
repeated this several times. Rather than treating the familiar ways to construct these
objects, he focused on the interplay between the separate structure components, 
i.e., order, topology, additive group, and multiplicative structure. For each basic structure,
he asked how they could be characterized, what happened when the characteristic 
properties were relaxed, what the homomorphic images were, and in which ways the 
structure components could be combined. Also related number systems such as $p$-adic 
numbers were treated in the same vein. H\"ahl and L\"owen worked out the lecture notes 
after these courses \cite{zahlen1}, \cite{zahlen2}, 
which became rather widespread at least within Germany. Many years 
later and after completion of the monograph \cite{CPP} on compact projective planes, 
he joined the same two, together with 
Grundh\"ofer, in producing an enhanced and expanded book version of those lecture notes 
\cite{classicalfields}.

In a more geometric vein, Salzmann wrote a book chapter on non-associative double 
loops and ternary fields together with Grundh\"ofer \cite{gr-sz}.

\section {Other types of geometries and Salzmann's research group}

Right from the start, Salzmann had a much broader vision of topological geometry, 
beyond the realm of projective planes. He suggested to his students areas of research where 
a fruitful program similar to his own could be successfully started. We already mentioned 
stable planes, where he did pioneering work before he left the area to his 
offspring of several generations (Strambach, Groh, L\"owen, Seidel, L\"owe, Stroppel). 
Circle geometries (M\"obius, Laguerre and Minkowski planes) became a very active area 
(Strambach, W\"olk, Schenkel, F\"ortsch, Steinke, Polster, Schroth). Generalized polygons 
and their relationship with circle geometries, and, more generally, topological buildings 
were treated by Grundh\"ofer, Joswig, Kramer, Schroth, Knarr, Van Maldeghem and other people 
not necessarily belonging to Salzmann's research group. Within the field 
of projective planes, some areas were mainly treated by group members, e.g., 
totally disconnected planes (Grundh\"ofer), translation planes (Betten, H\"ahl, L\"owe)
and the so-called shift planes (Betten, Knarr, L\"owen).
Differentiable planes were the last topic to be introduced (Otte, B\"odi, Immervoll, 
L\"owen, Pupeza).
Relationships to (not necessarily Riemannian) symmetric spaces were 
studied by L\"owen, Seidel, L\"owe, and
spatial geometries were treated by Groh, Betten, L\"owen.
Moreover, geometric objects occurring within planes and spaces were 
studied, such as ovals, unitals and parallelisms (Buchanan, H\"ahl, Stroppel, Betten, L\"owen). 
Strambach initiated a thorough 
investigation of groups of projectivities (compositions of central projections 
from one line to another).  This had been a major tool from the beginning, but now 
became a goal in its own right (Betten,  L\"owen, Grundh\"ofer, see \cite{loewindsh}).

\

There are many joint papers by Salzmann with his students, but mostly the latter 
worked freely in their own areas. The interaction became closer with his last group 
of students (Stroppel, Hubig, L\"uneburg, B\"odi, Priwitzer, Au\ss enhofer), 
compare Section \ref{higher}. 

It was a special feat of his strength that Salzmann left his students enough room 
so that they could develop their own substantial body of results. They kept close 
contact with him but were independent at the same time. A few times it happened that 
he himself got stuck in his program and some student removed the obstacles by introducing
completely new ideas. Nevertheless, in the area where his unique power lay, intricate 
exploitations of the interplay between group and geometry, 
keeping ideas straight even in the most involved situations,
piling up small steps until the
final result was reached --- in this domain he always stayed far ahead of his students,
right until the end of his life.

\section {Homogeneous planes}

Although there is still no proof that compact connected projective planes are manifolds, 
a theorem of Szenthe \cite{szenthe} implies that a locally compact group acting transitively on 
the (locally contractible) point set is a Lie group. 
Szenthe's proof is incomplete, but the result was proved 
later independently by Hofmann and Kramer \cite{hof-kra} and others. Montgomery's 
theorem \cite{mont} then shows that there is a compact transitive Lie group.
This allowed Salzmann to prove in \cite{sz75a} that 
a projective plane of dimension at most 16 admitting a point transitive group is classical. 
He used the limited topological information available at the time, compare the next 
section, and this accounts for the dimension restriction. His proof is based on a scrutiny 
of all feasible compact Lie groups. 
A completely different approach, based on a result of Hopf and Samelson 
\cite{hopf-sa} expressing the 
Euler characteristic of a homogeneous space in terms of Weyl groups,
allowed L\"owen to remove the dimension restriction \cite{loe81d}. 
In \cite{sz75b}, Salzmann obtained in addition that
an affine plane with a line transitive group is classical. In \cite{CPP}, Section 63, these 
results are presented, including  a comparison of the different approaches.
For stable planes with a group transitive on flags, a similar result was obtained by L\"owen 
\cite{loe83a}. Classical in this case means that the plane is the affine or projective
plane over $\BR$, $\BC$, $\BH$ or $\BO$, or the interior open orbit of the corresponding 
hyperbolic motion group, i.e., the open ball of unit radius in the affine plane. 
The method here was to use the reflections contained in the group in order 
to endow the point set with the structure of a symmetric space. Without recourse to the 
theory of symmetric spaces, the same result was proved by Salzmann in \cite{retrospect}.

\section {Higher dimensions}\label{higher}

The main questions that need to be solved before one can tackle the classification program
for planes with a group of large dimension
concern the topology of the plane (one hopes that lines are spheres of dimension 1,2,4, or  8 
like in the classical planes), the topology of the automorphism group (is it locally compact?),
the applicability of Lie theory, and stiffness. The difficulties arise from the fact that 
a higher dimensional ternary lacks internal structure. Looking at the classical cases,
one is lured into thinking that there should be a structure of real vector space or at 
least some weak substitute thereof. But there is no such thing. 

\

Almost the only  
general information is that ternaries are locally homogeneous and
locally contractible, as observed already 
by Freudenthal. Contractions are obtained using multiplication with the elements of a 
path that joins 1 to 0. This fact remained the foundation upon which everything was built.
As a consequence, one obtains that all relevant spaces are absolute neighborhood retracts and
Cantor manifolds (subsets of codimension $\ge 2$ cannot disconnect them), and closed 
subsets of the full dimension contain a non-empty open subset. In the strong form 
needed here, this was proved by Seidel \cite{seidel}.

\

Around 1980, after his second sabbatical 
with Dugundji in Los Angeles, and with his help, Salzmann knew a lot more about the 
first question, but for example he still could not rule out the possibility 
of 7-dimensional lines. 
Luckily, the first two questions were removed by his students. The group is always 
locally compact \cite{loe76b}, \cite{gr86} (in particular, its size 
can be measured by its topological dimension), 
and if the point set has finite dimension, 
then lines are homology manifolds of one of the classical dimensions \cite{loe83b}. 
So Salzmann was left with the really tough
questions of stiffness and the Lie property. 

\

Stiffness is measured by the stabilizer of a quadrangle, i.e., by the
automorphism group of a coordinatizing ternary $K$. 
In the classical cases, these are the groups $\mathop\mathrm {Aut} \BH = \mathrm{SO}_3$ and 
$\mathop\mathrm {Aut} \BO = \mathrm{G}_2$, the compact exceptional Lie group of dimension 14.
In general, it would help if one knew that $\mathop\mathrm {Aut} K$ is compact 
like in the classical cases, but this could be proved only in special cases.
The difficulties showed clearly in the article \cite{sz79b} 
on automorphisms of an 8-dimensional ternary $K$. 
The easiest case arises when there is no proper closed subternary: then 
there are no automorphisms of $K$, and hence the stabilizer 
of a quadrangle is trivial. If there happen to be subternaries of dimension 2 and 4, 
then something is known about how their stabilizers can act on them,  and the kernels of 
these actions have additional properties. However, there might be finite subternaries 
and none of dimension 1 or 2, and a lot of perseverance is needed in order to get 
the sharp bound  $\dim \mathop\mathrm {Aut} K \le 14$. If the bound is attained, then 
$\mathop\mathrm {Aut} K \cong \mathrm{G}_2$, but this does not imply that the associated 
plane is the octonion plane; there are examples due to Plaumann and Strambach 
which are not even translation planes 
\cite{p-s74}. Later B\"odi \cite{boe95}, \cite{boe94c} showed 
that $\dim \mathop\mathrm {Aut} K \le 11$
in all other cases, and there are many additional results for special situations, see  
\cite{update}, 2.6. 

\

The question about the Lie property of $\Sigma$ becomes 
somewhat less hard if one restricts attention 
to the connected component $\Delta$ of this group, where the task is to exclude the existence 
of small central subgroups of dimension zero. For 8-dimensional planes, Priwitzer 
\cite{pri94a} showed that $\Delta$ is a Lie group if its dimension is at least 12, a small 
bound if compared to the classical case, where $\dim \Sigma = 35$. In her proof, she considers
three types of collineations $\zeta$ in the center of $\Delta$ using distinct methods: 
those that have a center and an axis, those whose fixed points form a subplane, 
and the remaining ones, for which there is always a point $a$ such that $\zeta$ does not 
fix the line joining $a$ to its image $a^\zeta$. A similar result for 16-dimensional planes
is obtained in her joint paper with Salzmann \cite{pri-sz98}. In that case, 
$\dim \Delta \ge 27$ suffices (to be compared to 78, the dimension of the 
automorphism group of the octonion plane).

In the wider field of stable planes,
a different approach to the Lie group question was carried through successfully by Stroppel 
\cite{strop93e}, \cite{strop94}.
He used approximation of locally compact groups by Lie groups in order to make all 
the benefits of Lie theory available. However, Salzmann preferred to work with true Lie groups,
although he loved projective limits (and ultrapowers).

\

In pursuit of the classification program in dimension 8 and 16, tools from algebraic topology,
theory of Lie groups and their representations as well as transformation groups 
became increasingly important, and Salzmann used them skilfully. In particular, upper
bounds for the dimension of compact groups on a space of 
known dimension (Mann \cite{mann}) were 
much used. The classifications of doubly transitive transformation groups 
(Tits \cite{tits55}) and of 
compact transitive groups on spheres (Poncet \cite{poncet}) 
were equally useful, as well as Newman's 
theorem \cite{newm} (the fixed point sets of compact group actions have empty interior). 
A compilation of such helpful results  is given in the appendix of \cite{CPP}.

\

We already mentioned the technique of producing elations as products of reflections. 
This has been generalized and strengthened in several stages so as to become a powerful tool for
finding large (commutative) groups of elations with a common axis. One stage in this 
development is given by Salzmann in \cite{sz73a}, concerning planes of small dimension (2 or 4). 
The most powerful version is due to H\"ahl \cite{hh81}. It asserts among other things that 
a Lie group $\Delta$ fixing a line $A$ and containing a non-trivial homology with axis $A$ and 
center $c \notin A$ also contains an elation group with axis $A$ whose dimension 
equals that of the orbit $c^\Delta$.

A related sequence of results was obtained by 
Salzmann \cite{sz70}, \cite{sz85}, \cite{sz90}. The final version 
is the following. If $A$ is a line such that the elation groups  with axis $A$ and center 
$a \in A$ have equal dimensions except for one point $c \in A$, then the elation group
with axis $A$ and the special center $c$ is transitive (i.e., as large as possible). 
In particular, 
if the groups with center $a \in A$ have equal dimensions without exception, then the full
elation group with axis $A$ is transitive. 

\section {Translation planes}

Higher dimensional translation planes are the realm of H\"ahl. They
figure prominently in Salzmann's classification results, hence we shall give a few
examples that are both special and easily described. In general, translation planes 
are coordinatized by quasi-fields, resembling skew fields but lacking associativity 
of multiplication and one of the two distributive laws. 
The non-commutativity of quaternion 
and octonion multiplication allows to construct the \it mutations \rm $\BH_t$ and $\BO_t$
by retaining addition and introducing a new multiplication 
   $$x \circ_t y = t xy + (1-t) yx,$$
where $t \ne \frac 1 2$ is a real parameter. This yields a  semifield (= division algebra), 
i.e., a quasifield with both distributive laws. The choice of $t = 0$ yields the 
classical semifield, and the planes defined by the other mutations have 
automorphism groups  of dimension 17  and 40, respectively; moreover, the 
planes over the mutations $\BO_t$ are the only translation planes with a group of 
dimension at least 40. This is proved by H\"ahl \cite{hh88}, \cite{hh90},  
see also \cite{CPP} 82.23 and 82.27. In fact, L\"uneburg \cite{luen92} showed that
the last statement is true without the restriction to translation planes. 

\

The only 8-dimensional non-classical translation planes admitting a group with the 
largest possible dimension 18 are obtained from the following semifields: 
Take a complex number $h$ from the first quadrant of the unit circle and with 
positive real part. Split every quaternion uniquely as $a = r(a) + p(a)\cdot h$, 
where $r(a)\in \BR$ and $p(a)$ is a  pure quaternion. Then set
    $$a \circ_h x = r(a) \cdot x + p(a) \cdot x \cdot h.$$ 
The translation planes with a group of dimension 18 are the planes over these 
semifields, their dual planes, and the planes obtained from them by a technique 
called transposition. This was proved by H\"ahl \cite{hh75c}, \cite{hh86a}; a proof 
is also given in \cite{CPP} 82.2.   

\

Nearfields are quasifields with associative multiplication. According to 
Tits \cite{tits52}, \cite{tits56},
the only locally compact, connected proper nearfields are those constructed by 
Kalscheuer \cite{ka40}. They are obtained from the quaternions by 
introducing the new multiplication 
    $$a \circ x = a\cdot \vert a\vert^{-ir}\cdot x \cdot \vert a \vert ^{ir}$$
whith $r >0$ a real parameter. The corresponding planes have groups of dimension 17, see
\cite{CPP} 82.24. We remark that Tits showed that, apart from the classical 
affine linear groups on $\BR$, $\BC$, $\BH$, these nearfields define the only 
sharply 2-transitive Lie transformation groups on manifolds.

\section {Hughes planes} \label{Hughes}

If $\BF$ is one of the classical skew fields and $\BD$ is the classical division algebra 
of twice the dimension, then the special linear group $\mathop {\mathrm {SL}_3}\BF$ 
acts on the projective plane over $\BD$, leaving the Baer subplane defined by $\BF$ invariant. 
We mentioned Salzmann's result that the plane over $\BC$ is entirely determined by this 
action (in fact, by this group itself) if $\BF = \BR$. 
In \cite{sz81a} and \cite{sz82a}, he showed that there 
is a one-parameter family of \it Hughes planes \rm admitting this action in each of the 
higher dimensions, and that there are no other possibilities in dimension 8. The name was 
given to the planes because Hughes had constructed very similar examples in the finite case. 
The 8-dimensional Hughes planes were known as abstract 
projective planes from Dembowski's paper \cite{demb71a},
and they can be obtained by extending the idea of homogeneous coordinates to coordinates 
taken from a Kalscheuer nearfield, see Maier and Stroppel \cite{m-str}. 
The 16-dimensional ones have to 
be constructed entirely from the group.  H\"ahl \cite{hh86b} showed that they, too, 
are characterized by the group. Salzmann's sketchy proof of their existence and 
their topological properties was completed by Niemann \cite{niemann}. See also Section 86 of 
\cite{CPP}, where a full account of topological Hughes planes is given. The automorphism groups 
of the non-classical Hughes planes have dimension 17 and 36 for $\BF = \BC$ and 
$\BF = \BH$, respectively.

\section {State of the art}

Planes of dimension 8 with a group of dimension at least 17 are either classical or 
Hughes planes or translation planes (up to duality). 
This result of Salzmann \cite{sz79a}, \cite{sz81a}
is proved in \cite{CPP} 84.28, using simplifications due to Stroppel \cite{strop93e}. 

\

After the completion of the monograph \cite{CPP}, Salzmann published 16 more research 
papers in print. At the age of 83 and 86, he wrote two articles 
giving a summary of the achievements in the 
8- and 16-dimensional cases, respectively, including also new 
results and proofs \cite{retrospect}, \cite{update}. 
In the 8-dimensional case, for every combination of
9 possible fixed element configurations and 4 possible structural features of a 
group (semi-simple, normal torus, normal vector group, general), it is listed 
which lower bounds on 
the dimension $d$ of the group imply that the plane is known or is classical or is a Hughes 
plane etc.. For example, if there are no fixed elements and the group is semi-simple, 
then dimension $d \ge 12$ implies that the plane is a Hughes plane.  Without semi-simplicity, 
all planes with $d \ge 12$ are known, and  $d \ge 18$ implies classical. In general, 
planes with $d \ge 17$ are known and $d \ge 19$ implies classical. 

\

Among planes of dimension 16, we already mentioned that 
the most homogeneous non-classical ones are 
the planes over mutations of the octonions, with groups of dimension 40. 
In \cite{CPP}, Section 87, 
the path leading to this result is indicated. 
Contributors include Salzmann, Stroppel, B\"odi, Hubig,
Priwitzer, L\"uneburg. A summary of results like in the 8-dimensional case is given in 
\cite{update}. As a surprise, in one of two joint papers by 
H\"ahl and Salzmann \cite{hh-sz10} there 
appeared a class of 16-dimensional planes over cartesian fields
and with a group of dimension 38; they are not translation planes and are reminiscent 
of the 2-dimensional Cartesian planes encountered right at the beginning.

\section {Books and Surveys}

A survey on stable and projective planes, mainly of dimension 2, 
was given by Salzmann \cite{sz67b}.
The known results on groups of projectivities of topological projective and circle planes 
were reviewed by L\"owen in the proceedings of a conference on von 
Staudt's point of view in Geometry
\cite{loewindsh}. In the same proceedings, Salzmann reported on what could be learned 
from projectivities about the topology of projective planes \cite{sz81b}.
Grundh\"ofer and Salzmann \cite{gr-sz} contributed a report on topological double loops 
(generalizing linear ternaries).  A survey on 
linear topological geometries in general, including spatial geometries and relationships 
with symmetric spaces is \cite{handb}, complemented by a similar report on circle 
geometries by Steinke \cite{steinke}.
Also there is a book by  Polster and Steinke on linear and 
circle geometries on surfaces \cite{polstein}. The most comprehensive account of compact 
projective planes is the monograph \cite{CPP}. 
The book about the classical fields \cite{classicalfields} 
has also been mentioned before. Finally, there is Salzmann's survey on Baer subplanes 
both in abstract and in topological projective planes \cite{sz03a}.

\bibliographystyle{plain}

\begin{thebibliography}{9}

\bibitem{adams}
Adams, J. F.,
On the non-existence of elements of Hopf invariant one, 
Ann. of Math. 72, 20--104, 1960.

\bibitem{arens}
Arens, R.,
Topologies for homeomorphism groups, 
Amer. J. Math. 68, 593--610, 1946.

\bibitem{betten70}
Betten, D.,
Nicht-desarguessche $4$-dimensionale Ebenen, 
Arch. Math. 21, 100--102, 1970.
 

\bibitem{betten72b}
Betten, D.,
$4$-dimensionale Translationsebenen, 
Math. Z. 128, 129--151, 1972.  

\bibitem{betten73a}
Betten, D.,
$4$-dimensionale Translationsebenen mit $8$-dimensionaler 
Kollineationsgruppe, 
Geom. Dedicata 2, 327--339, 1973.
         
\bibitem{betten73b}
Betten, D.,         
$4$-dimensionale Translationsebenen mit irreduzibler 
Kollineationsgruppe, 
Arch. Math. 24, 552--560, 1973.
         
\bibitem{betten73c}
Betten, D.,      
Die komplex-hyperbolische Ebene, 
Math. Z. 132, 249--259, 1973.

\bibitem{boe95}
B\"odi, R.,
On the dimensions of automorphism groups of 8-dimensional ternary fields I,
J. Geom. 52, No. 1-2, 30-40, 1995. 

\bibitem{boe94c}
B\"odi, R.,
On the dimensions of automorphism groups of eight-dimensional ternary fields II,
Geom. Dedicata 53, No. 2, 201-216, 1994. 

\bibitem{brouwer09}
Brouwer, L. E. J.,
Die Theorie der endlichen kontinuierlichen Gruppen, unabh\"angig von
den Axiomen von Lie, 
Math. Ann. 67, 246--267, 1909.
         
         
\bibitem{brouwer12}
Brouwer, L. E. J.,         
Beweis des ebenen Translationssatzes, 
Math. Ann. 72, 37--54, 1912.

\bibitem{demb71a}
Dembowski, P.,
Generalized Hughes planes, 
Canad. J. Math. 23, 481--494, 1971.


\bibitem{fr57a}
Freudenthal, H.,
Kompakte projektive Ebenen, 
Illinois J. Math. 1, 9--13, 1957.
	        
\bibitem{fr57b}
Freudenthal, H.,	        
Zur Geschichte der Grundlagen der Geometrie. Zugleich eine Besprechung der
8. Aufl. von Hilberts \lq\lq Grundlagen der Geometrie\rq\rq, 
Nieuw Arch. Wisk. (4) 5, 105--142, 1957.

\bibitem{gr86}
Grundh\"ofer, T.,
Automorphism groups of compact projective planes, 
Geom. Dedicata 21, 291--298 1986.




\bibitem{handb}
T. Grundh\"ofer, T. and  L\"owen, R.,
Linear topological geometries,
in: F. Buekenhout (ed.), 
Handbook of Incidence Geometry, Chapter 23, pp. 1255 -- 1324,
Amsterdam: North Holland, 1995.

\bibitem{hh81}
H\"ahl, H.,
Homologies and elations in compact, connected projective planes, 
Topology Appl. 12, 49--63, 1981.

\bibitem{hh86b}
H\"ahl, H.,
Charakterisierung der kompakten, zusammenh\"angenden 
Moufang-Hughes-Ebenen anhand ihrer Kollineationen, 
Math. Z. 191, 117--136, 1986.

\bibitem{hh88}
H\"ahl, H.,
Die Oktavenebene als Translationsebene mit gro\ss{}er
Kollineationsgruppe, 
Monatsh. Math. 106, 265--299, 1988.

\bibitem{hh90}
H\"ahl, H.,
Sechzehndimensionale lokalkompakte Translationsebenen, deren
Kollineationsgruppe $\mathrm G_2$ enth\"alt, 
Geom. Dedicata 36, 181--197, 1990.

\bibitem{hh75c}
H\"ahl, H.,
Geometrisch homogene vierdimensionale reelle Divisionsalgebren, 
Geom. Dedicata 4, 333--361, 1975.

\bibitem{hh86a}
H\"ahl, H.,
Achtdimensionale lokalkompakte Translationsebenen mit mindestens
$17$-dimensionaler Kollineationsgruppe, 
Geom. Dedicata 21, 299--340, 1986.

\bibitem{hilb}
Hilbert, D.,
Grundlagen der Geometrie, 
Leipzig: Teubner, 1899.
7th edition:
Leipzig und Berlin: Teubner, 1930.

\bibitem{hof-kra}
Hofmann, K.H. and Kramer, L.,
Transitive actions of locally compact groups on locally contractible spaces, 
J. Reine Angew. Math. 702, 227-243, 2015; erratum ibid. 702, 245-246, 2015. 

\bibitem{hopf-sa}
Hopf, H. and Samelson, H.,
Ein Satz \"uber die Wirkungsr\"aume geschlossener Liescher Gruppen,
Comment. Math. Helv. 13, 240--251, 1941.

\bibitem{ka40}
Kalscheuer, F.,
Die Bestimmung aller stetigen Fastk\"orper \"uber dem K\"orper der
reellen Zahlen als Grundk\"orper, 
Abh. Math. Sem. Hansische Univ. (Hamburg) 13, 413--435, 1940.

\bibitem{loe76b}
L\"owen, R.,
Vierdimensionale stabile Ebenen, 
Geom. Dedicata 5, 239--294, 1976.

\bibitem{loewindsh}
L\"owen, R.,
Projectivities and the geometric structure of topological planes, 
in: P. Plaumann, K. Strambach (eds.), 
Geometry --- von Staudt's point of view, 
Proc. Bad Windsheim 1980, pp. 339--372,
Dordrecht etc.: Reidel, 1981.

\bibitem{loe81d}
L\"owen, R.,
Homogeneous compact projective planes, 
J. Reine Angew. Math. 321, 217--220, 1981.

\bibitem{loe83a}
L\"owen, R.,
Stable planes with isotropic points, 
Math. Z. 182, 49--61, 1983.

\bibitem{loe83b}
L\"owen, R., 
Topology and dimension of stable planes: On a conjecture of H. 
Freudenthal, 
J. Reine Angew. Math. 343, 108--122, 1983.

\bibitem{luen92}
L\"uneburg, M.,
Involutionen, aufl\"osbare Gruppen und die 
Klassifikation topologischer Ebenen, 
Mitt. Math. Sem. Gie\ss{}en 209, 1992.


\bibitem{kol}
Kolmogoroff, A.,
Zur Begr\"undung der projektiven Geometrie, 
Ann. of Math. 33, 175--176, 1932.

\bibitem{m-str}
Maier, P. and Stroppel, M.,
Pseudo-homogeneous coordinates for Hughes planes,
Can. Math. Bull. 39, No. 3, 330-345, 1996. 

\bibitem{mann}
Mann, L. N.,
Dimensions of compact transformation groups, 
Michigan Math. J. 14, 433--444, 1967.


\bibitem{mont}
Montgomery, D.,
Simply connected homogeneous spaces, 
Proc. Amer. Math. Soc. 1, 467--469, 1950.


\bibitem{mozip}
Montgomery, D. and Zippin, L.,
Topological transformation groups, 
New York: Interscience, 1955. 

\bibitem{newm}
Newman, M. H. A.,
A theorem on periodic transformations of spaces, 
Quart. J. Math. Oxford 2, 1--8, 1931.

\bibitem{niemann}
Niemann, K.,
Geometrie und Topologie der sechzehndimensionalen Moufang-Hughes-Ebenen,
Staatsexamensarbeit, Universit\"{a}t Kiel, 1990. 

\bibitem{pick}
Pickert, G.,
Projektive Ebenen, 
Berlin etc.: Springer, 1955.

\bibitem{p-s70}
Plaumann, P. and Strambach, K., 
Zusammenh\"angende Quasik\"orper mit Zentrum, 
Arch. Math. 21,455--465, 1970.  

\bibitem{p-s74}
Plaumann, P. and Strambach, K.,
Hurwitzsche Tern\"ark\"orper, 
Arch. Math. 25, 129--134, 1974.

\bibitem{poncet}
Poncet, J.,
Groupes de Lie compacts de transformations de l'espace euclidien 
et les sph\`eres comme espaces homog\`enes, 
Comment. Math. Helv. 33, 109--120, 1959.

\bibitem{polstein}
Polster, B. and Steinke, G.,
Geometries on surfaces,
Encyclopedia of Mathematics and Its Applications. 84. Cambridge: Cambridge 
University Press, 2001. 

\bibitem{pont}
Pontrjagin, L. S.,
\"Uber stetige algebraische K\"orper, 
Ann. of Math. 33, 163--174, 1932.

\bibitem{pri94a}
Priwitzer, B.,
Large automorphism groups of 8-dimensional projective planes are Lie groups,
Geom. Dedicata 52, No. 1, 33-40, 1994. 




\bibitem{rich}
Richardson, R. W.,
Groups acting on the $4$-sphere, 
Illinois J. Math. 5, 474--485, 1961.

\bibitem{sz-z}
Salzmann, H. and Zeller K.,
Singularit\"aten unendlich oft differenzierbarer Funktionen,
Math. Z. 62, 354-367, 1955. 

\bibitem{sz55}                                                    
Salzmann, H., 
\"Uber den Zusammenhang in topologischen projektiven Ebenen, 
Math. Z. 61, 489--494, 1955.
         
\bibitem{sz57}                                                       
Salzmann, H., 
Topologische projektive Ebenen, 
Math. Z. 67, 436--466, 1957.
	 
\bibitem{sz58}
Salzmann, H., 
Kompakte zweidimensionale projektive Ebenen, 
Arch. Math. 9, 447--454, 1958.
	 
\bibitem{sz59a}
Salzmann, H., 
Homomorphismen topologischer projektiver Ebenen, 
Arch. Math. 10, 51--55, 1959.

\bibitem{sz59b}
Salzmann, H., 	 
Topologische Struktur zweidimensionaler projektiver Ebenen, 
Math. Z. 71, 408--413, 1959.
         
\bibitem{sz59c}
Salzmann, H.,          
Viereckstransitivit\"at der kleinen projektiven Gruppe einer Moufang-Ebene, 
Illinois J. Math. 3, 174--181, 1959.

\bibitem{sz62a}                                                 
Salzmann, H.,         
Kompakte zweidimensionale projektive Ebenen, 
Math. Ann. 145, 401--428, 1962.

\bibitem{szutrecht}
Salzmann, H.,
Topologische projektive Ebenen,
Algebr. Topol. Foundations Geom., Proc. Colloq. Utrecht, August 1959, 157-163, 1962. 

\bibitem{sz62b}                                                 
Salzmann, H., 	 
Kompakte Ebenen mit einfacher Kollineationsgruppe, 
Arch. Math. 13, 98--109, 1962.

\bibitem{sz63a}                                                  
Salzmann, H., 	 
Characterization of the three classical plane geometries, 
Illinois J. Math. 7, 543--547, 1963.

\bibitem{sz63b}                                                   
Salzmann, H.,          
Zur Klassifikation topologischer Ebenen, 
Math. Ann. 150, 226--241, 1963.
	 
\bibitem{sz64}                                                    
Salzmann, H., 	 
Zur Klassifikation topologischer Ebenen. II, 
Abh. Math. Sem. Univ. Hamburg 27, 145--166, 1964.   
            
\bibitem{sz65}                                                      
Salzmann, H., 
Zur Klassifikation topologischer Ebenen. III, 
Abh. Math. Sem. Univ. Hamburg 28, 250--261, 1965.

\bibitem{bu-sz}                                                         
Busemann, H., Salzmann, H.,
Metric collineations and inverse problems,
Math. Z. 87, 214-240, 1965. 


\bibitem{sz66}
Salzmann, H.,  	 
Polarit\"aten von Moulton-Ebenen, 
Abh. Math. Sem. Univ. Hamburg 29, 212--216, 1966.

\bibitem{sz67a}
Salzmann, H., 	 
Kollineationsgruppen ebener Geometrien, 
Math. Z. 99, 1--15, 1967.

\bibitem{sz67b}                                                     
Salzmann, H., 	 
Topological Planes, 
Adv. Math. 2, 1--60, 1967.

\bibitem{szcollill}
Salzmann, H.,
Topological geometries. 
Proc. Project. Geom. Conf. Univ. Illinois 1967, 119-120, 1967. 

\bibitem{sz69a}
Salzmann, H., 	 
Geometries on surfaces, 
Pacific J. Math. 29, 397--402, 1969.

\bibitem{sz69b}                                                  
Salzmann, H., 	 
Kompakte vier-dimensionale Ebenen, 
Arch. Math. 20, 551--555, 1969.

\bibitem{sz69c}
Salzmann, H., 	 
Homomorphismen komplexer Tern\"ark\"orper, 
Math. Z. 112, 23--25, 1969.

\bibitem{sz70}
Salzmann, H., 	                                                     
Kollineationsgruppen kompakter, vier-dimensionaler Ebenen, 
Math. Z. 117, 112--124, 1970.

\bibitem{sz71}
Salzmann, H., 	 
Kollineationsgruppen kompakter $4$-dimensionaler Ebenen. II,          
Math. Z. 121, 104--110, 1971.

\bibitem{zahlen1}
Salzmann, H.,
Zahlbereiche. Teil I: Die reellen Zahlen. Vorlesung,
Universit\"at T\"ubingen, 1971. 

\bibitem{sz72a}
Salzmann, H., 	 
$4$-dimensional projective planes of Lenz type III, 
Geom. Dedicata 1, 18--20, 1972.

\bibitem{sz72b}
Salzmann, H., 	                                                   
Homogene $4$-dimensionale affine Ebenen, 
Math. Ann. 196, 320--322, 1972.

\bibitem{sz72c}
Salzmann, H., 	 
Baer-Unterebenen $4$-dimensionaler Ebenen, 
Arch. Math. 23, 337--341, 1972.

\bibitem{sz73a}
Salzmann, H., 	 
Elations in four-dimensional planes, 
Gen. Top. Appl. (=Topology Appl.) 3, 121--124, 1973.

\bibitem{sz73b}
Salzmann, H., 	                                                 
Kompakte, vier-dimensionale projektive Ebenen mit
$8$-dimensionaler Kollineationsgruppe, 
Math. Z. 130, 235--247, 1973.

\bibitem{zahlen2}
Salzmann, H.,
Zahlbereiche. Teil II: Die rationalen Zahlen. Teil III: Die komplexen Zahlen. Vorlesung,
Mathematisches Institut der Universit\"at T\"ubingen, 1973. 

\bibitem{sz73c}
Salzmann, H.,                                                          
Reelle Kollineationen der komplexen projektiven Ebene, 
Geom. Dedicata 1, 344--348, 1973.

\bibitem{inmemdembowski}
Salzmann, H.,
In memoriam: Peter Dembowski. 
Proc. internat. Conf. projective Planes, Washington State Univ. 1973, 1-6, 1973. 

\bibitem{sz74}
Salzmann, H., 	 
Compact planes of Lenz type III, 
Geom. Dedicata 3, 399--403, 1974.

\bibitem{sz75a}
Salzmann, H., 	                                                
Homogene kompakte projektive Ebenen, 
Pacific J. Math. 60, 217--234, 1975.   

\bibitem{sz75b}
Salzmann, H.,                                                                               
Homogene affine Ebenen, 
Abh. Math. Sem. Univ. Hamburg 43, 216--220, 1975.

\bibitem{szbwg}
Salzmann, H.,
Zur nicht-euklidischen Geometrie,
Abh. Braunschw. wiss. Ges. 27, 119, 1977. 

\bibitem{sz79a}
Salzmann, H., 	 
Compact $8$-dimensional projective planes with large collineation      
groups, 
Geom. Dedicata 8, 139--161, 1979.

\bibitem{sz79b}
Salzmann, H.,                                                                   
Automorphismengruppen $8$-dimensionaler Tern\"ark\"orper, 
Math. Z. 166, 265--275, 1979.

\bibitem{sz81a}                                                         
Salzmann, H., 	 
Kompakte, $8$-dimensionale projektive Ebenen mit gro\ss{}er
Kollineationsgruppe, 
Math. Z. 176, 345--357, 1981.

\bibitem{sz81b}                                                             
Salzmann, H.,          
Projectivities and the topology of lines, 
in: P.~Plaumann, K.~Strambach (eds.), 
Geometry --- von Staudt's point of view, 
Proc. Bad Windsheim 1980, pp. 313--337, 
Dordrecht etc.: Reidel, 1981.

\bibitem{sz82a}                                                         
Salzmann, H.,          
Baer-Kollineationsgruppen der klassischen projektiven Ebenen, 
Arch. Math. 38, 374--377, 1982.

\bibitem{sz82b}
Salzmann, H.,          
Compact $16$-dimensional projective planes with large collineation 
groups, 
Math. Ann. 261, 447--454, 1982.

\bibitem{loe-sz}
L\"owen, R. and Salzmann, H.,
Collineation groups of compact connected projective planes,
Arch. Math. 38, 368-373, 1982. 


\bibitem{sz83}
Salzmann, H.,          
Compact $16$-dimensional projective planes with large collineation 
groups. II, 
Monatsh. Math. 95, 311--319, 1983.

\bibitem{sz84}
Salzmann, H.,          
Compact $16$-dimensional projective planes with large collineation 
groups. III, 
Math. Z. 185, 185--190, 1984.

\bibitem{sz85}                                                      
Salzmann, H.,          
Homogeneous translation groups, 
Arch. Math. 44, 95--96, 1985.

\bibitem{sz87}
Salzmann, H., 	 
Compact $16$-dimensional projective planes with large collineation 
groups IV, 
Canad. J. Math. 39, 908--919, 1987.

\bibitem{sz90}                                                       
Salzmann, H.,          
Compact $8$-dimensional projective planes, 
Forum Math. 2, 15--34, 1990.

\bibitem{gr-sz}
Grundh\"ofer, T. and Salzmann, H.
Locally compact double loops and ternary fields, 
Chap. XI of: O.~Chein, H.~Pflugfelder, J.~Smith (eds.), 
Quasigroups and Loops: Theory and Applications, pp. 313--355, 
Berlin: Heldermann, 1990. 


\bibitem{CPP}
H. Salzmann, D. Betten, T. Grundh\"ofer, H. H\"ahl, R. L\"owen, and M. Stroppel,
Compact projective planes, 
Walter de Gruyter, Berlin, New York, 1995.

\bibitem{sz98}
Salzmann, H.,
Characterization of 16-dimensional Hughes planes, 
Arch. Math. 71, No. 3, 249-256, 1998. 

\bibitem{pri-sz98}
Priwitzer, B. and Salzmann, H.,
Large automorphism groups of 16-dimensional planes are Lie groups. 
J. Lie Theory 8, No. 1, 83-93, 1998. 

\bibitem{sz99a}
Salzmann, H.,
Compact 16-dimensional projective planes,
Result. Math. 35, No. 1-2, 192-196, 1999. 

\bibitem{sz99b}
Salzmann, H.,
Large automorphism groups of 16-dimensional planes are Lie groups II,
J. Lie Theory 9, No. 2, 481-486, 1999. 

\bibitem{sz00}
Salzmann, H.,
On the classification of 16-dimensional planes,
Beitr. Algebra Geom. 41, No. 2, 557-568, 2000. 

\bibitem{sz01}
Salzmann, H.,
Near-homogeneous 16-dimensional planes,
Adv. Geom. 1, No. 2, 145-155, 2001. 


\bibitem{sz03a}
Salzmann, H.,
Baer subplanes,
Ill. J. Math. 47, No. 1-2, 485-513, 2003. 


\bibitem{sz03b}
Salzmann, H.,
16-dimensional compact projective planes with 3 fixed points,
Adv. Geom. 2003, Spec. Issue, S153-S157, 2003. 

\bibitem{sz05}
Salzmann, H.,
16-dimensional compact projective planes with a large group fixing two points and two lines, 
Arch. Math. 85, No. 1, 89-100, 2005. 

\bibitem{hh-sz05}
H\"ahl, H. and Salzmann, H.,
16-dimensional compact projective planes with a large group fixing two points and two lines, 
Arch. Math. 85, No. 1, 89-100,  2005. 


\bibitem{classicalfields}
Salzmann, H., Grundh\"ofer, T., H\"ahl, H., L\"owen, R.,
The classical fields. Structural features of the real and rational numbers. 
Encyclopedia of Mathematics and Its Applications 112. Cambridge: Cambridge University Press,
2007.

\bibitem{sz08}
Salzmann, H.,
16-dimensional compact projective planes with a collineation group of dimension $\ge 35$,
Arch. Math. 90, No. 3, 284-288, 2008. 


\bibitem{sz10}
Salzmann, H.,
Classification of 8-dimensional compact projective planes,
J. Lie Theory 20, No. 4, 689-708, 2010. 

\bibitem{hh-sz10}
H\"ahl, H. and Salzmann, H.,
16-dimensional compact projective planes with a large group fixing two 
points and only one line, 
Innov. Incidence Geom. 11, 213-235, 2010. 

\bibitem{sz14}
Salzmann, H.,
8-dimensional compact planes with an automorphism group which has a normal vector subgroup, 
J. Lie Theory 24, No. 1, 123-146, 2014. 

\bibitem{sz16a}
Salzmann, H.,
Reminiszenzen an G\"unter Pickert,
J. Geom. 107, No. 2, 221-224, 2016.


\bibitem{sz16b}
Salzmann, H.,
Semi-simple groups of compact 16-dimensional planes,
J. Geom. 107, No. 2, 249-255, 2016. 

\bibitem{sz19}
Salzmann, H.,
Groups of compact 8-dimensional planes: conditions implying the Lie property,
Innov. Incidence Geom. 17, No. 3, 201-220, 2019. 

\bibitem{retrospect}
Salzmann, H., 
Compact planes, mostly 8-dimensional. A retrospect, arxiv:1402.0304, 2014.

\bibitem{update}
Salzmann, H., 
Compact 16-dimensional planes. An update, arxiv:1706.03696, 2017

\bibitem{seidel}
Seidel, H.-P.,
Locally homogeneous ANR-spaces, 
Arch. Math. 44, 79--81, 1985.


\bibitem{sk54}
Skornjakov, L. A.,
Topological projective planes (Russian), 
Trudy Moskov Mat. Obshch. 3, 347--373, 1954. 

\bibitem{sk57}
Skornjakov, L. A.,
Systems of curves on a plane (Russian), 
Trudy Moskov Mat. Obshch. 6, 135--164, 1957. 

\bibitem{steinke}
Steinke, G. F.,
Topological circle geometries,
in: F. Buekenhout (ed.), 
Handbook of Incidence Geometry, Chapter 24, pp. 1325 -- 13524,
Amsterdam: North Holland, 1995.

\bibitem{strop93e}
Stroppel, M.,
Stable planes with large groups of automorphisms: the interplay of incidence,
topology, and homogeneity, 
Habilitationsschrift, Darmstadt 1993.

\bibitem{strop94}
Stroppel, M.,
Lie theory for non-Lie groups,
J. Lie Theory 4, No. 2, 257-284, 1994. 

\bibitem{szenthe}
Szenthe, J.,
On the topological characterization of transitive Lie group actions, 
Acta Sci. Math. (Szeged) 36, 323--344, 1974.

\bibitem{tits52}
Tits, J.,
Sur les groupes doublement transitifs continus, 
Comment. Math. Helv. 26, 203--224, 1952.

\bibitem{tits55}
Tits, J.,
Sur certaines classes d'espaces homog\`enes de groupes de Lie, 
M\'em. de l'Acad\'emie Royale de Belgique, 
Classe de Sciences XXIX, Fasc. 3, 1955.

   
\bibitem{tits56}
Tits, J.,
Sur les groupes doublement transitifs continus: correction et
compl\'ements, 
Comment. Math. Helv. 30, 234--240, 1956.

\end{thebibliography}

\noindent{Rainer L\"owen\\ 
Institut f\"ur Analysis und Algebra\\
Technische Universit\"at Braunschweig\\
Universit\"atsplatz 2\\
38106 Braunschweig\\
Germany

\end{document}